\documentclass[12pt]{article}
\usepackage[latin1]{inputenc}
\usepackage[english]{babel}
\usepackage{latexsym}
\usepackage{amsmath}
\usepackage{amssymb}
\usepackage{amsfonts}
\usepackage{a4wide}
\newtheorem{thm}{Theorem}
\title{\bf How many quasiplatonic surfaces?}
\author{Jan--Christoph Schlage--Puchta,  J\"urgen Wolfart\\
{\small Mathematisches Institut der Universität, Eckerst. 1, D--79104 Freiburg}\\
{\small Mathematisches Seminar der Universit\"at,  Postfach 11 19 23, 
D--60054 Frankfurt a.M.}\\
{\small \tt jcp@math.uni-freiburg.de \quad
  wolfart@math.uni-frankfurt.de}}
\date{}
\begin{document}
\maketitle

\begin{abstract}
We show that the number of isomorphism classes of quasiplatonic Riemann
surfaces of genus $\,\le g\,$ has a growth of type $\,g^{\log g}\,$. The
number of non--isomorphic regular dessins
of genus $\,\le g\,$ has the same growth type.\\
MSC Index: 20E07, 30F10, 14H30.
\end{abstract}

Quasiplatonic Riemann surfaces $X$ of genus $\,g>1\,$ can be characterized
in many different ways (see e.g. [Wo2, Thm. 4], or [Si2]),  
e.g. by the property that the orders of their
automorphism groups are isolated local maxima in the corresponding
moduli space of Riemann surfaces of genus $g\,$. For the present paper 
we will only use the equivalent statement that their
universal covering groups $\Gamma$ are torsion free normal subgroups of 
finite index in some Fuchsian triangle groups $\,\Delta= \Delta(p,q,r)\,$ of signatures 
$\,(p,q,r)\,$, and that conversely the quotient $\,\Gamma \backslash {\bf H}\,$ of
the upper half plane ${\bf H}$ by any finite index 
torsion free normal subgroup $\Gamma$ of a Fuchsian 
triangle group is a quasiplatonic surface $X\,$. The dessin it carries
can be described as the quotient by $\Gamma$ of a $\Delta$--invariant
tesselation of ${\bf H}\,$, with $\,\Delta/\Gamma\,$ acting as a group
of automorphisms on $X$ and on the dessin. 

Notations: Let $\,R(g;p,q,r)\,$ be the number of normal torsion free
subgroups $\Gamma$ of genus $g$ in $\,\Delta= \Delta(p,q,r)\,$, and let 
\[ R(g) \; := \; \sum_{p,q,r} R(g;p,q,r) \; , \quad S(g) \; := \;
\sum_{1< \gamma \le g} R(\gamma) \]
be the number of all non--isomorphic regular dessins of this genus $g$ and its summatory
function. We will show first 

\begin{thm}
There are constants $\,g_0,c_1,c_2>0\,$ such that for all genera $\,g>g_0\,$
\[ g^{c_1 \log g} \; < \; S(g) \; < \; g^{c_2 \log g} \;. \]
\end{thm} 

To prove the lower bound, take three different primes
$\,p,q,r\,$ giving a Fuchsian triangle group $\,\Delta= \Delta(p,q,r)\,$ with presentation 
\[ \Delta \;=\; \langle \,\gamma_0 \,,\, \gamma_1 \,| \, \gamma_0^p =
\gamma_1^q = (\gamma_0 \gamma_1)^r = 1 \, \rangle \;. \]
All normal subgroups of index $\,n>1\,$ are torsion free. By the 
Riemann--Hurwitz formula, their genus $g$ is related to $n$ via 
\[ 84(g-1) \; \ge \; n \; = \; (2g-2) \, (1-\frac{1}{p}-\frac{1}{q}-\frac{1}{r})^{-1} \; >
\; 2g-2 \;, \] 
and by [MSP] we have a lower bound for the summatory growth function 
\[  \sum_{1< 2(\gamma-1)<n} R(\gamma;p,q,r) \; \ge 
 \; | \{ \Gamma \lhd \Delta(p,q,r) \, \mid \,1< (\Delta(p,q,r) : \Gamma) \le n \,\}| \; >
 \;  n^{c_1\log n} \] 
for all $\,n>n_0\,$ for some $n_0$ depending on $\,p,q,r\,$. 
Taking only that term in the sum $\, R(\gamma) =
\sum_{p,q,r} R(\gamma;p,q,r)\,$ coming from the triangle group 
$\,\Delta= \Delta(p,q,r)\,$ for the chosen prime triple signature we obtain 
\[ S(g) \; \ge \; \sum_{1< \gamma \le g} R(\gamma;p,q,r) \;  \ge \;
 | \{ \Gamma \lhd \Delta \, \mid \, 1<(\Delta :\Gamma) \le 2(g-1) \,\}|  \]
\[ \qquad \ge \;  | \{ \Gamma \lhd \Delta \, \mid \, 1<(\Delta :\Gamma) \le g
\,\}| \; > \; g^{c_1 \log g}  \] 
for all $\,g \ge n_0\,$. 

To prove the upper bound, recall Lubotzky's estimate $\, \nu^{6 (\Omega(\nu)+1)}\,$
for the number of index $\nu$ normal subgroups in the free
group with two generators ([LS, Thm. 2.7]) where $\Omega(\nu)$ 
denotes the number of prime divisors of $\nu$ counted with multiplicities.
For any fixed Fuchsian triangle group $\,\Delta= \Delta(p,q,r)\,$ it  
implies
\[  | \{ \Gamma \lhd \Delta \, \mid \, (\Delta : \Gamma) \le n \,\}| \; \le \;
\sum_{\nu \le n}  \nu^{6 (\Omega(\nu)+1)} \; < \; n^{c_3\log n}  \] 
for some constant $c_3\,$, since $\,\Omega(\nu) \le \log_2 \nu\,$. 
This upper bound is {\em a fortiori} valid 
for the torsion free normal subgroups, hence
\[ \sum_{1<\gamma \le g} R(\gamma;p,q,r) \; < \; (84g)^{c_3 \log (84g)} \;.\] 
Since $\,\Delta/\Gamma\,$ has generators of orders $\,p,q,r\,$, we have
moreover $\,R(\gamma;p,q,r) = 0\,$ for $\,p,q\,$ or $\,r > 84g\,$ ($\,>
|\Delta/\Gamma|\,$), therefore 
\[ S(g) \; = \; \sum_{1 < \gamma \le g} R(\gamma) \; < \; \sum_{p,q,r \le
  84g} (84g)^{c_3 \log (84g)} \; < \; (84g)^{3+c_3 \log (84g)} \; < \;
g^{c_2 \log g}  \]
for all $\,g>g_0\,$ with suitable $c_2$ and $g_0\,$. 

\begin{thm}
Let $\,Q(g)\,$ denote the number of isomorphism classes of quasiplatonic
Riemann surfaces of genera $\,\gamma, \; 1<\gamma \le g\,$. 
With the same constants $\,g_0,c_1,c_2>0\,$ 
as in Theorem 1 we have for all $\,g> g_0\,$ 
\[ \frac{1}{120} \, g^{c_1 \log g} \; < \; Q(g) \; < \; g^{c_2 \log g} \;. \]
\end{thm}

{\em Proof.} Since every quasiplatonic surface is obtained from a 
regular dessin (equivalently, from a torsion free normal subgroup in a
triangle group), and is uniquely determined by that dessin, the upper 
bound follows from Theorem 1. The lower bound
follows similarly, but a quasiplatonic surface can be obtained by 
up to five different types of regular dessins ([Gi]), and
another overcount can happen: in a fixed triangle group $\,\Delta=
\Delta(p,q,r)\,$ several torsion free normal subgroups can be 
${\rm PSL}_2({\bf R})$--conjugate, leading to isomorphic surfaces. 
In [GW], Thms. 5, 6, 7, it is shown
that such conjugations take place in a finite extension of
$\Delta$ which is again a triangle group. By Singerman's work [Si1] 
the maximal possible index
between Fuchsian triangle groups is known to be $24$ (occurring for $\, \Delta(2,3,7) \supset
\Delta(7,7,7) \,$), so 
we have at most $120$ normal subgroups counted in the proof of Theorem 1 
leading to isomorphic surfaces (by
a more detailed analysis, this number can considerably decreased). This
gives the lower bound for $\,Q(g)\,$. \\
Another consequence of Thm. 1 is  

\begin{thm}
With the same constants $\,g_0,c_1,c_2>0\,$ as in Theorem 1, 
the number of non--isomorphic regular dessins in genus $\,g> g_0\,$ is
$\,R(g) < g^{c_2 \log g}\;$. Infinitely often we have
\[  g^{-1+c_1 \log g} \; < \; R(g) \;.\]
An analogous statement holds for the number of quasiplatonic
surfaces of genus $g\,$.
\end{thm}

{\em Remarks.} 1) From the tables in [Wo1, Sec. 6] of regular dessins in genera $\,g \le
4\,$ and work of Kuribayashi and Kimura [KK] for $\,g=5\,$ one may
deduce 
\[ S(5) = 104 \qquad \mbox{and} \qquad Q(5) = 37 \;,\]
to be compared with $ \, 5^{\log 5} \approx 13\,$. [GAP] calculations
[Sc] indicate that also for $\,5 < g \le 10\,$ one has always 
\[ g^{\log g} < S(g) < g^{2 \log g} \;.\]
2) Counting regular dessins in genera $0$ and $1$ is different from
higher genera. In genus $0$ there is only one surface but an infinity of
regular dessins defined by the trivial triangle groups of signatures $\,
(1,n,n), \, (2,2,n)\,$ and those corresponding to the platonic bodies,
i.e. $\,(2,3,3),\, (2,3,4),\,(2,3,5)\,$. \\
In genus $1$ the torsion free normal subgroups of the triangle groups of
signatures\\ 
$\,(3,3,3), \,(2,3,6),\, (2,4,4)\,$ give infinitely many
non--isomorphic elliptic curves with regular dessins, but all fall in
two isogeny classes only, see [SiS].

{\em Acknowledgement.} The authors thank the referee for several useful remarks.

\subsection*{References}

[GAP] The GAP Group: GAP --- Groups, Algorithms, and Programming,\\ 
{\tt  http://www.gap-system.org} 

[Gi] E. Girondo: Multiply quasiplatonic Riemann surfaces, Experimental
Mathematics 12 (4) (2003), 463--475 

[GW] E. Girondo, J. Wolfart: Conjugators of Fuchsian groups and
quasiplatonic surfaces, preprint 2004,  {\tt
  http://www.math.uni-frankfurt.de/$\sim$wolfart}, to appear in
Quarterly J. Math. Oxford 

[KK] A. Kuribayashi, H. Kimura: Automorphism Groups of Compact Riemann
Surfaces of Genus Five, Journal of Algebra 134 (1990), 80--103 

[LS] A. Lubotzky, D. Segal: Subgroup growth, PM 212, Birkh\"auser 2003

[MSP] T.W. M\"uller, J.--Chr. Schlage--Puchta: Normal growth of large
groups, II, to appear in Archiv der Mathematik 

[Sc] M. Schmitt: Regul\"are Hypermaps, Diplomarbeit, Frankfurt 2001

[Si1] D. Singerman: Finitely Maximal Fuchsian Groups, J. London
Math. Soc. (2) 6 (1972), 29--38 

[Si2] D. Singerman: Riemann Surfaces, Belyi Functions and Hypermaps, 
pp. 43--68 in Topics on Riemann Surfaces and Fuchsian Groups (ed:
E. Bujalance, A.F. Costa, E. Martinez), London Math. Soc. LNS 287,
Cambridge UP 2001

[SiS] D. Singerman, R.I. Syddall: Belyi Uniformization of Elliptic
Curves, Bull. London Math. Soc. 139 (1997), 443--451

[Wo1] J. Wolfart: Triangle groups and Jacobians of CM type, Frankfurt
2000,\\ 
{\tt  http://www.math.uni-frankfurt.de/$\sim$wolfart} 

[Wo2] J. Wolfart: ABC for polynomials, dessins d'enfants, and
uniformization --- a survey, to appear in the proceedings of the ELAZ 
conference 2004, ed.: W. Schwarz, J. Steuding, 
{\tt  http://www.math.uni-frankfurt.de/$\sim$wolfart} 

\end{document}